\documentclass[12pt]{article}

\newcommand{\mathC}{{\mathbb C}}

\newcommand{\bbX}{{\bf X}}

\newcommand{\bbY}{{\bf Y}}

\newcommand{\bbu}{{\bf u}}
\newcommand{\bbv}{{\bf v}}

\usepackage{graphics}
\usepackage{amsmath,amssymb}
\usepackage[dvips]{xy}
\xyoption{all}

\begin{document}
\baselineskip=6.0mm

\newcommand{\ignore}[1]{}{}

\renewcommand{\theequation}{\arabic{section}.\arabic{equation}}

\newcommand{\lbl}{\label}


\newcommand{\eq}[1]{$(\ref{#1})$}

\newcommand{\al}{\alpha}                         
\newcommand{\bt}{\beta}                          
\newcommand{\ga}{\gamma}                         
\newcommand{\Ga}{\Gamma}                         
\newcommand{\de}{\delta}                         
\newcommand{\De}{\Delta}                         
\newcommand{\ep}{\epsilon}                       
\newcommand{\ve}{\varepsilon}                    
\newcommand{\la}{\lambda}                        
\newcommand{\La}{\Lambda}                        
\newcommand{\ro}{\rho}                           
\newcommand{\ta}{\tau}                           
\newcommand{\si}{\sigma}                         

\newcommand{\be}{\begin{equation}}               
\newcommand{\ee}{\end{equation}}                 
\newcommand{\bea}{\begin{eqnarray}}              
\newcommand{\eea}{\end{eqnarray}}                
\newcommand{\bean}{\begin{eqnarray*}}            
\newcommand{\eean}{\end{eqnarray*}}              
\newcommand{\ba}{\begin{array}}                  
\newcommand{\ea}{\end{array}}                    
\newcommand{\nn}{\nonumber}                      
\newcommand{\mb}{\mbox}                          

\newcommand{\f}{\frac}

\newcommand{\ra}{\rightarrow}                    
\newcommand{\llra}{\longleftrightarrow}          

\newcommand{\stac}{\stackrel}                    
\newcommand{\noin}{\noindent}                    

\newcommand{\qed}{\nobreak\quad\vrule width6pt depth3pt height10pt}

\newcommand{\ngi}{n \ra \infty}

\pagestyle{myheadings} \markright{Entropy method for the left tail}

\thispagestyle{plain}

\begin{center}
{\Large\bf Entropy method for the left tail}
\vskip 0.2in

Hyungsu Kim,\footnote{Department of Mathematics, Yonsei
University, Seoul, Korea 120-749.

\hspace{0.1cm} E-mail: gudtn@yonsei.ac.kr}$^{,a}$ Chul Ki
Ko,\footnote{University College, Yonsei University, Seoul, Korea
120-749. E-mail: kochulki@yonsei.ac.kr } and Sungchul
Lee\footnote{Department of Mathematics, Yonsei University, Seoul,
Korea 120-749.

\hspace{0.1cm}E-mail: sungchul@yonsei.ac.kr.

\hspace{0.1cm}$^{a}$ Supported by the BK21 project of the
Department of Mathematics, Yonsei University.}$^{,a}$ \vskip 0.2in
\vskip 0.2in
\end{center}

\begin{abstract}
When we use the entropy method to get the tail bounds, typically
the left tail bounds are not good comparing with the right ones.
Up to now this asymmetry has been observed many times.
Surprisingly we find an entropy method for the left tail that
works in the resembling way that it works for the right tail. This
new method does not work in all the cases. We provide a meaningful
example.

\ignore{
\vspace{1.0cm}
\hrule
\vspace{0.1cm}
\noindent {\bf \it AMS $1991$ subject classification.}
                           Primary 60D05;
                           secondary 05C80, 90C27.
\newline
\noindent {\bf \it  Key words and phrases  } Rate of convergence,
minimal matching,
 minimal spanning tree,  traveling salesman problem.
\vspace{0.1cm}
\hrule
}
\end{abstract}

\section{Introduction.}\setcounter{equation}{0}

\quad In recent years, interesting developments took place in the
analysis of the spectrum of large random matrices. In particular,
the asymptotic distribution of the largest eigenvalue has been a
subject of hot interest.

Let $\bbX=(X_{ij})$ be an $n\times n$ complex hermitian matrix such
that the entries $X_{ij}$ on and above diagonal are independent
complex (real on the diagonal) centered normal random variables with
variance $1$. Let $\la_1\ge\la_2\ge\cdots\ge\la_n$ be the $n$ real
eigenvalues of $\f{1}{\sqrt{n}}\bbX$. There have been many
researches of the concentration of the largest eigenvalue $\la_1$ or
the concentration of the $k$-th largest eigenvalue $\la_k$.
Regarding the concentration of the $k$-th largest eigenvalue, we
know of three results; Alon, Krivelevich, Vu (2002), Meckes (2004),
Maurer (2006). Alon, Krivelevich, Vu (2002) and Meckes (2004) used
Talagrand's method whereas Maurer (2006) used the entropy method.
Since our main theme of this paper is the  entropy method, we state
Maurer's concentration result.\\

\noin {\bf Theorem.}[Maurer (2006)]  Let $\bbX=(X_{ij})$ be an
$n\times n$ real symmetric matrix such that the entries $X_{ij}$
on and above diagonal are independent with $|X_{ij}|\le 1$. Let
$\la_1\ge\la_2\ge\cdots\ge\la_n$ be the $n$ real eigenvalues of
$\bbX$. Then, for all $k$, $n\ge 1$, and for all $t\ge 0$, \be
P(\la_k-E\la_k\ge t)\le \exp\left(-\f{t^2}{16k^2}\right),\
P(\la_k-E\la_k\le -t)\le
\exp\left(-\f{t^2}{16k^2+2kt}\right).\lbl{1.5} \ee

The left tail bounds in \eq{1.5} are larger than the right ones.
This asymmetry usually happens when we use the entropy method to get
the tail bounds. However, this asymmetry is not observed in the
works of Alon, Krivelevich, Vu (2002) and Meckes (2004) which are
based on Talagrand's method. In these works the left tail bounds are
same to the right ones. In addition, the centering is the median not
the mean. This symmetry and the centering are typical with
Talagrand's method.

In this paper we found an entropy method for the left tail that
works in the resembling way that it works for the right tail by
controlling the term $\De^2$ carefully (see \eq{2.9} in Section 2
for the definition of $\De^2$), and  give a meaningful example.

The rest of the paper is organized as follows. In Section 2, we
develop an entropy method for the left tail. In Section 3, we apply
this new method to the interesting case including the $k$-th largest
eigenvalue.
\section{Entropy method for the left tail.} \setcounter{equation}{0}

\quad The concentration of measure phenomenon for the product
measures has been investigated in depth by Talagrand (1995, 1996)
in a most remarkable way. His method has been applied to various
interesting cases. In many cases his method made new-record
concentration inequalities and in some cases his method even
produced non-trivial concentration inequalities for the first
time. However, his method is technically too complicated. Hence
many people tried to simplify his proof and studied to find an
alternative to reproduce  and more ambitiously to extend his
result. One of the successful alternatives is the entropy method.
Here we explain the minimum details of the entropy method to show
our contributions on this interesting subject. See Ledoux (1996),
Massart (2000), Boucheron, Lugosi, Massart (2000, 2003), Maurer
(2006) for the full details.

Let $X_1,\ldots,X_n$ be independent and let $G=G(X_1,\ldots,X_n)>0$.
Define the entropy $H(G)$ and the partial entropy $H_k(G)$ by
\bea
H(G)&:=&EG\log G-EG\log EG,\nn\\
H_k(G)&:=&E_kG\log G-E_kG\log E_kG,\nn
\eea
where $E$ is the integration over $X_1,\ldots,X_n$
whereas $E_k$ is the integration over $X_k$ only.
So, the entropy $H(G)$ is a real number
but the partial entropy $H_k(G)$ is a random variable which does not depends on $X_k$.

Some classical formulas of the entropy are quite helpful;
\bea
H(G)&=&\sup_{T}EG(\log T-\log ET),\lbl{2.1}\\
H(G)&=&\inf_{c}EG(\log G-\log c)-(G-c),\lbl{2.2}
\eea
where the
supremum in \eq{2.1} is taken over the strictly positive random
variables $T$ and where the infimum in \eq{2.2}  runs over the
strictly positive constants $c$. \eq{2.1} is called the duality
formula of the entropy and \eq{2.2} is called the variation
formula.

Here is the well-known entropy inequality (or tensorization inequality)
which follows from the duality formula \eq{2.1}.\\

\noin
{\bf Lemma 1.} [Entropy inequality]
\be
H(G)\le\sum_{k=1}^nEH_k(G).\lbl{2.3}
\ee

Now, let $Z=Z(X_1,\ldots,X_n)$ be the random variable of interest.
We apply the entropy inequality to the random variable $e^{\la Z}$.
Then, we have
\be
E\la Ze^{\la Z}-Ee^{\la Z}\log Ee^{\la Z}\le\sum_{k=1}^nEH_k(e^{\la Z}).\lbl{2.4}
\ee
To estimate the term $EH_k(e^{\la Z})$, we apply the variation formula \eq{2.2} to the partial entropy $H_k(e^{\la Z})$;
$
H_k(e^{\la Z})=\inf_{c}E_ke^{\la Z}(\la Z-\log c)-(e^{\la Z}-c).
$
Since the integration $E_k$ is only over $X_k$,
during the evaluation of the partial entropy $H_k(e^{\la Z})$
we can treat
all the other random variables $X_j$, $1\le j\neq k\le n$, as fixed constants.
So, in fact $c$ can be chosen as a function of $X_1,\ldots,X_{k-1},X_{k+1},\ldots,X_n$,
or even as a function of $X_1,\ldots,X_{k-1},X_k',X_{k+1},\ldots,X_n$,
where $X_k'$ is an independent copy of $X_k$ and $X_k'$ is independent to $X_1,\ldots,X_n$.
This subtle point on $c$ is crucial for the further development of the theory.
If we choose a particular ``constant'' $c_0$ to estimate the partial entropy $H_k(e^{\la Z})$,
then we have
\be
H_k(e^{\la Z})\le E_ke^{\la Z}(\la Z-\log c_0)-(e^{\la Z}-c_0).\lbl{2.5}
\ee
To get a good concentration inequality,
we have to choose $c_0$ well-designed for the random variable $Z$ of interest.

There are many possible choices of $c_0$.
Massart (2000) and Boucheron, Lugosi, Massart (2003) chose
$$
c_0:=\exp\left(\la Z(X_1,\ldots,X_k',\ldots,X_n)\right):=e^{\la Z_k},
$$
where $X_k'$ is an independent copy of $X_k$.

Boucheron, Lugosi, Massart (2000) chose
$$
c_0:=\exp\left(\la Z(X_1,\ldots,\hat X_k,\ldots,X_n)\right):=e^{\la Z_k}.
$$
Here, $\hat X_k$ means that we drop out $X_k$ from the argument of $Z$.
In other word, we evaluate the value of $Z$ based not $\{X_1,\ldots,X_n\}$ but $\{X_1,\ldots,X_n\}\setminus\{X_k\}$.
This is possible because of the special nature of the random variable $Z$ they considered.

Maurer (2006) chose
\be
c_0:=\exp\left(\la \inf_{x_k}Z(X_1,\ldots,x_k,\ldots,X_n)\right):=e^{\la Z_k},\lbl{2.6}
\ee
where the infimum runs over all the possible values $x_k$
which $X_k$ can take as a function value or
over a compact set containing the support of the distribution of $X_k$.
He used this $c_0$ (or $Z_k$) to get the right tail bound in Theorem A.
He also use the same $Z_k$ to obtain the left tail bound in the same Theorem.

In this paper we follow the footsteps of Maurer for the right tail
bound. However, to get a better left tail bound we choose the
following $c_0=e^{\la Z_k}$ for the left tail bound;
\be
c_0:=\exp\left(\la\sup_{x_k}Z(X_1,\ldots,x_k,\ldots,X_n)\right)
:=e^{\la Z_k},\lbl{2.7}
\ee
where the supremum runs over a compact
set containing the support of the distribution of $X_k$. This
choice does not always come with a sensible $\De^2$ (see \eq{2.9}
below for the definition of $\De^2$). However, in many cases with
this choice we do have $\De^2$ with $\|\De^2\|_{\infty}<\infty$.

Let's recall what we have done so far with the entropy inequality.
We first apply the entropy inequality to $G=e^{\la Z}$
where $Z$ is the random variable of interest.
Then, the term $EH_k(e^{\la Z})$ appears in the inequality.
To estimate the term $EH_k(e^{\la Z})$,
with a particular choice $c_0=e^{\la Z_k}$ we apply the variation formula to $H_k(e^{\la Z})$.
Then, we get the following log-Sobolev inequality.\\

\noin
{\bf Lemma 2.} [Log-Sobolev inequality]
If $-\la (Z-Z_k)\le 0$ for all $k$, then
\be
E\la Ze^{\la Z}-Ee^{\la Z}\log Ee^{\la Z}\le\f{\la^2}{2}Ee^{\la Z}\De^2,\lbl{2.8}
\ee
where
\be
\De^2:=\sum_{k=1}^n\left(Z-Z_k\right)^2.\lbl{2.9}
\ee

{\bf Proof.}
With a particular choice $c_0=e^{\la Z_k}$, from \eq{2.5} we have
\bea
H_k(e^{\la Z})
&\le&E_ke^{\la Z}\left(e^{-\la (Z-Z_k)}-\left(1-\la (Z-Z_k)\right)\right)\nn\\
&=&E_ke^{\la Z}\f{e^{-\la (Z-Z_k)}-\left(1-\la (Z-Z_k)\right)}{\la^2(Z-Z_k)^2}\la^2(Z-Z_k)^2.\nn
\eea
If $-\la (Z-Z_k)\le 0$ for all $k$,
since $(e^x-(1+x))/x^2$ is an increasing function with the function value $1/2$ at the trouble spot $x=0$,
and (hence) since $(e^x-(1+x))/x^2\le 1/2$ for $x\le 0$, we have then
$$
H_k(e^{\la Z})\le\f{\la^2}{2}E_ke^{\la Z}(Z-Z_k)^2.
$$
Plug this estimate into \eq{2.4} and we get the log-Sobolev inequality \eq{2.8}.
\qed\\

To distinguish our choice \eq{2.7} from Maurer's choice \eq{2.6}, from now on we let
\bea
\De_M^2&:=&\sum_{k=1}^n\left(Z-\inf_{x_k}Z(X_1,\ldots,x_k,\ldots,X_n)\right)^2:=\sum_{k=1}^n\left(Z-Z_k^{(M)}\right)^2,\nn\\
\De_L^2&:=&\sum_{k=1}^n\left(Z-\sup_{x_k}Z(X_1,\ldots,x_k,\ldots,X_n)\right)^2:=\sum_{k=1}^n\left(Z-Z_k^{(L)}\right)^2.\nn
\eea

Here is our entropy method for the left tail,
which is a simple consequence of the log-Sobolev inequality.\\

\noin {\bf Theorem 1.} (i) If $\|\De_M^2\|_{\infty}\le\infty$,
then for $t\ge 0$ \be P(Z-EZ\ge t)\le
\exp\left(-\f{t^2}{2\|\De_M^2\|_{\infty}}\right).\lbl{2.10} \ee

(ii) If $\|\De_L^2\|_{\infty}\le\infty$, then for $t\ge 0$
\be
P(Z-EZ\le -t)\le\exp\left(-\f{t^2}{2\|\De_L^2\|_{\infty}}\right).\lbl{2.11}
\ee

{\bf Remark.} As Maurer pointed out in private communication,
$\|\De_M^2\|_{\infty}\neq \|\De_L^2\|_{\infty}$. However, in
practice we don't know the exact values of $\|\De_M^2\|_{\infty}$
and $\|\De_L^2\|_{\infty}$. Instead we calculate the upper bounds
of $\|\De_M^2\|_{\infty}$ and $\|\De_L^2\|_{\infty}$. In case
$\|\De_M^2\|_{\infty}=\|\De_L^2\|_{\infty}<\infty$,  \eq{2.10} and
\eq{2.11} provide the same left and right tail bounds.

{\bf Proof.} The right tail bound \eq{2.10} is Theorem 1 of Maurer
(2006). So, we can safely skip its proof. In fact, the left tail
bound \eq{2.11} also follows from the same argument, the so-called
Herbst's argument. For reader's convenience here we reproduce the
Herbst's argument to get \eq{2.11}.

In this proof, we will use only the negative $\la\le 0$.
Then, (since by our choice of $Z_k^{(L)}$, $Z-Z_k^{(L)}\le 0$)
we have  $-\la (Z-Z_k)\le 0$ for all $k$.
So, we can use the log-Sobolev inequality \eq{2.8}.
Since $\|\De_L^2\|_{\infty}\le\infty$, by \eq{2.8}
$$
E\la Ze^{\la Z}-Ee^{\la Z}\log Ee^{\la Z}\le\f{\la^2}{2}\|\De_L^2\|_{\infty}Ee^{\la Z}.
$$

Divide the both sides by $\la^2Ee^{\la Z}$.
Then, we have
$$
\f{d}{d\la}\f{1}{\la}\log Ee^{\la(Z-EZ)}\le\f{\|\De_L^2\|_{\infty}}{2}.
$$

Recall $\la\le 0$.
So, we integrate the both sides from $\la$ to 0.
Since $\la^{-1}\log Ee^{\la(Z-EZ)}\ra 0$ as $\la\ra 0$, we have then
$
-\la^{-1}\log Ee^{\la(Z-EZ)}\le-\|\De_L^2\|_{\infty}\la/2
$
or
\be
Ee^{\la(Z-EZ)}\le\exp\left(\f{\|\De_L^2\|_{\infty}}{2}\la^2\right).\lbl{2.12}
\ee

Now, by Chebyshev's inequality with the choice $\la=-t/\|\De_L^2\|_{\infty}\le 0$ we have
the left tail bound \eq{2.11}; by \eq{2.12},
$$
P(Z-EZ\le -t)
\le e^{\la t}Ee^{\la(Z-EZ)}
\le \exp\left(\la t+\f{\|\De_L^2\|_{\infty}}{2}\la^2\right)
=\exp\left(-\f{t^2}{2\|\De_L^2\|_{\infty}}\right).
$$
\qed

\section{Example.}\setcounter{equation}{0}

\quad In this section, we apply the entropy method for the left
tail (Theorem 1) to the eigenvalues of sample covariance matrix.
In a near future we hope to see  many more exciting examples.

Let $\bbX=(X_{ij})$ be an $n\times N$ complex matrix with the
independent entries $X_{ij}$. Let $\la_1\ge\la_2\ge\cdots\ge\la_n$
be the $n$ positive eigenvalues of $\f{1}{N}\bbX\bbX^*$. Then,
under the suitable condition on the distribution of $X_{ij}$ the
Mar\v{c}enko-Pastur theorem (Mar\v{c}enko and Pastur (1967)) says
that as $n\ra\infty$, $N\ra\infty$, $n/N\ra c (0<c<\infty)$, the
empirical spectral distribution
$\frac{1}{n}\sum_{k=1}^n\de_{\la_k}$ of the sample covariance
matrix $\f{1}{N}\bbX\bbX^*$ converges to the Mar\v{c}enko-Pastur
law. This time we use the Mar\v{c}enko-Pastur scaling. For the
sample covariance matrix we don't know any established
concentration inequality to compare with. So, it is rather natural
to work with the Mar\v{c}enko-Pastur scaling.
Here is our result.\\

\noin {\bf Theorem 2.} Let $\bbX=(X_{ij})$ be an $n\times N$
complex matrix with the independent entries $X_{ij}$, which are
bounded by 1, i.e., $|X_{ij}|\le 1$. Let
$\la_1\ge\la_2\ge\cdots\ge\la_n$ be the $n$ positive eigenvalues
of $\f{1}{N}\bbX\bbX^*$. Then, for all $k$, $n, N\ge 1$, and for
all $t\ge 0$, \be P(\la_k-E\la_k\ge t)\le
\exp\left(-\f{Nt^2}{2n^2}\right),\ P(\la_k-E\la_k\le -t)\le
\exp\left(-\f{Nt^2}{2n^2}\right).\nn \ee

{\bf Proof.}
Let $X_t$ be the $t$-th column of $\bbX$.
To denote the dependency of the eigenvalues on the matrix $\bbX$,
we let $\la_1(\bbX)\ge\la_2(\bbX)\ge\cdots\ge\la_n(\bbX)$ be the $n$ positive eigenvalues of $\f{1}{N}\bbX\bbX^*$.
Fix $1\le k\le n$ and let $Z:=Z(\bbX):=\la_k(\bbX)$
be the $k$-th largest eigenvalue of $\f{1}{N}\bbX\bbX^*$.

Fix $1\le t_0\le N$. From the given $n\times N$ matrix $\bbX$
delete the $t_0$-th column $X_{t_0}$ and add $x_{t_0}$ where
$x_{t_0}$ is a constant column vector of size $n$ whose entries
are all bounded by 1. Call this new $n\times N$ matrix as $\bbY$.
Using this $\bbY$ we define $Z^{(M)}_{t_0}$ by
\be
Z^{(M)}_{t_0}:=\inf_{\bbY}Z(\bbY)=\inf_{x_{t_0}}Z(\bbY).\lbl{3.7}
\ee

Let $S^{k}$ be an arbitrary $k$-dimensional complex linear
subspace of $\mathbb{C}^{n}$. By the Courant-Fischer
representation theorem (look up Theorem 7.7 of Zhang (1999) for
the Courant-Fischer representation theorem),
\bea Z(\bbX)
&=&\f{1}{N}\max_{S^{k}}\min_{\bbv \in S^{k}, \bbv^*\bbv=1}\bbv^*\bbX\bbX^*\bbv\nn\\
&=&\f{1}{N}\max_{S^{k}}\min_{\bbv \in S^{k}, \bbv^*\bbv=1}\bbv^*\left(\sum_{t=1}^{N}X_tX_t^*\right)\bbv\nn\\
&=&\f{1}{N}\max_{S^{k}}\min_{\bbv \in S^{k}, \bbv^*\bbv=1}\bbv^*\left(\sum_{t=1}^{N}Y_tY_t^*+X_{t_0}X^*_{t_0}-x_{t_0}x^*_{t_0}\right)\bbv\nn\\
&=&\f{1}{N}\max_{S^{k}}\min_{\bbv \in S^{k}, \bbv^*\bbv=1}\bbv^*\left(\sum_{t=1}^{N}Y_tY_t^*\right)\bbv+\bbv^*\left(X_{t_0}X^*_{t_0}-x_{t_0}x^*_{t_0}\right)\bbv\nn\\
&\le&\f{1}{N}\max_{S^{k}}\min_{\bbv \in S^{k}, \bbv^*\bbv=1}\bbv^*\left(\sum_{t=1}^{N}Y_tY_t^*\right)\bbv+\max_{\bbu \in \mathC^n, \bbu^*\bbu=1}\bbu^*\left(X_{t_0}X^*_{t_0}-x_{t_0}x^*_{t_0}\right)\bbu\nn\\
&=&Z(\bbY)+\f{1}{N}\max_{\bbv \in \mathC^n,\bbv^*\bbv=1}\bbv^*\left(X_{t_0}X^*_{t_0}-x_{t_0}x^*_{t_0}\right)\bbv.\nn
\eea
Since $|X_{lt_0}|\le 1$ and since $\bbv^*\bbv=1$, we have
\bea
Z(\bbX)-Z(\bbY)
&\le&\f{1}{N}\max_{\bbv\in\mathC^n,\bbv^*\bbv=1}\bbv^*\left(X_{t_0}X^*_{t_0}-x_{t_0}x^*_{t_0}\right)\bbv\nn\\
&=&\f{1}{N}\max_{\bbv\in\mathC^n,\bbv^*\bbv=1}\left(\overline{X^*_{t_0}\bbv}\right)\left(X^*_{t_0}\bbv\right)-\left(\overline{x^*_{t_0}\bbv}\right)\left(x^*_{t_0}\bbv\right)\nn\\
&\le&\f{1}{N}\max_{\bbv\in\mathC^n,\bbv^*\bbv=1}\left|\sum_{l=1}^n\overline{X}_{lt_0}v_l\right|^2\nn\\
&\le&\f{1}{N}\max_{\bbv\in\mathC^n,\bbv^*\bbv=1}\left(\sum_{l=1}^n\left|X_{lt_0}\right|^2\right)\left(\sum_{l=1}^n\left|v_l\right|^2\right)\nn\\
&\le&\f{n}{N}.\nn
\eea
Take the infimum over $x_{t_0}$. Then, by the choice of $Z^{(M)}_{t_0}$ given in \eq{3.7} we have
$$
0\le Z-Z_{t_0}\le\f{n}{N}.
$$
So,
\be
\De_M^2:=\sum_{t_0=1}^N\left(Z-Z_{t_0}\right)^2\le\f{n^2}{N}.\lbl{3.8}
\ee
By \eq{3.8} and by Theorem 1 (i) we have the right tail bound
for $Z=\la_k$.

Now, we consider the left tail. When we choose $Z_{t_0}$, instead
of taking the infimum this time we take the supremum. Define
$Z^{(L)}_{t_0}$ by
\be
Z^{(L)}_{t_0}:=\sup_{\bbY}Z(\bbY)=\sup_{x_{t_0}}Z(\bbY).\lbl{3.9}
\ee
Then, by the Courant-Fischer representation theorem we have
$$
Z(\bbY)-Z(\bbX)\le\f{n}{N}.
$$
Take the supremum over $x_{t_0}$. Then, by the choice of
$Z^{(L)}_{t_0}$ given in \eq{3.9} we have
\be
0 \le
Z^{(L)}_{t_0}-Z\le\f{n}{N}.\nn
\ee
So,
\be
\De_L^2\le\f{n^2}{N}.\lbl{3.10}
\ee
By \eq{3.10} and by Theorem 1
(ii), we have the left tail bound for $Z=\la_k$. \qed

\vspace{1.0cm} \noin {\bf Acknowledgment.}
We appreciate Maurer's comment.
He pointed out
$\|\De_M^2\|_{\infty}\neq \|\De_L^2\|_{\infty}$ to us.

\end{document}